\documentclass{article}
\usepackage[utf8]{inputenc}
\usepackage{appendix}
\title{
Prawitz's Conjecture is False; So What?}
\author{Will Stafford}

\usepackage{setspace}
\usepackage[margin=1in]{geometry}
\usepackage[style=authoryear]{biblatex}
\bibliography{bib.bib}

\usepackage{xspace}
\usepackage{todonotes}
\usepackage{enumitem}

\usepackage{amsfonts}
\usepackage{amsmath}
\usepackage{amssymb}
\usepackage{amsthm}
\theoremstyle{definition}
\newtheorem{thm}{Theorem}

\newtheorem{defn}[thm]{Definition}

\newtheorem{ex}[thm]{Example}

\newtheoremstyle{named}{}{}{}{}{\bfseries}{.}{.5em}{\thmnote{#3}}
\theoremstyle{named}
\newtheorem*{nmd}{Theorem}

\newcommand{\set}[1]{\{#1\}}

\usepackage{ebproof}

\newcommand{\edit}[1]{#1}
\begin{document}
\maketitle
\begin{abstract}
    Several recent results bring into focus the superintuitionistic nature of most notions of proof-theoretic validity, but little work has been done evaluating the consequences of these results.  Proof-theoretic validity claims to offer a formal explication of how inferences follow from the definitions of logic connectives (which are defined by their introduction rules).  This paper explores whether the new results undermine this claim.  It is argued that, while the formal results are worrying, superintuitionistic inferences are valid because the treatments of atomic formulas are insufficiently general, and a resolution to this issue is proposed.
    
\end{abstract}
\tableofcontents
\section{Introduction}
%
%

Proof-theoretic validity was proposed by Prawitz as an explication of Gentzen's famous observation that the elimination rules of intuitionistic logic follow from the introduction rules. \edit{Proof-theoretic validity is a property that potential proofs have when they can, in some sense, be reduce to introduction rules alone. }Our goal in this paper is to assess whether recent work has shown that the formal explication of this idea as formulated by Prawitz is inadequate.  Let us start, though with the initial idea and Gentzen's famous quote:
\begin{quote}
    The introductions represent, as it were, the `definitions' of the symbols concerned, and the eliminations are no more, in the final analysis, than the consequences of these definitions. This fact may be expressed as follows: In eliminating a symbol, we may use the formula with whose terminal symbol we are dealing only `in the sense afforded it by the introduction of that symbol.'\phantom{aaaaaaaaaaaaaaaaaaaaaaaaaaa} \phantom{a}\hfill \hbox{(\cite[p. 189]{Gentzen1935} translated \cite[p. 295]{Gentzen1964})}
\end{quote}
Gentzen is suggesting that introduction rules for logical connectives should be treated as definitions, while (intuitionistic) elimination rules should be treated as consequences of the introduction rules.  In the quote above, however, `consequence' cannot mean logical consequence in the traditional model-theoretic sense. Rather, an introduction rule for a connective defines the connective by the premises, which are sufficient for proving a formula with the connective as the main connective, and the elimination rule for the connective is a consequence of the  introduction rule because it treats formulas with the connective as the main connective as the claim that premises of the  introduction rule can be proven.  

An illustration can help precisify intuitions.  Below are rules for the introduction and elimination of implication:

\begin{center}
    \begin{prooftree}
    \hypo{[A]}
    \infer[no rule]1{\vdots}
    \infer[no rule]1{B}
    \infer1{A\rightarrow B}
    \end{prooftree}\qquad\qquad
    \begin{prooftree}
    \hypo{A\rightarrow B}
    \hypo{A}
    \infer2{B}
    \end{prooftree}
\end{center}
When we examine the introduction rule for implication, we see that we can introduce an implication when we have a proof that takes $A$ as an assumption and transforms it into $B$.  This is a method for going from $A$ to $B$.  When we turn to the elimination rule, we see that it allows us to infer from `$A$ implying $B$' and $A$ holding that $B$ holds.  So, the elimination rule treats $A\rightarrow B$ in `the manner afforded by the introduction' because the introduction rule treated it as a \emph{method} for going from $A$ to $B$, and the elimination rule \emph{uses} it to go from $A$ to $B$.  

To formally explicate this idea, Prawitz connected Gentzen's idea with the notion of a detour and normalisation. A detour in a proof occurs when a connective is introduced and then eliminated.  For $\rightarrow$, a detour is depicted below. This detour can be removed via a reduction, which is also displayed.
\begin{center}
    \begin{prooftree}
    \hypo{[A]}
    \infer[no rule]1{\vdots}
    \infer[no rule]1{B}
    \infer1{A\rightarrow B}
    \hypo{\vdots}
    \infer[no rule]1{A}
    \infer2{B}
    \end{prooftree}\qquad$\Rightarrow$\qquad
    \begin{prooftree}
    \hypo{\vdots}
    \infer[no rule]1{A}
    \infer[no rule]1{\vdots}
    \infer[no rule]1{B}
    \end{prooftree}
\end{center}
Note that this reduction symbolically captures the relationship we discussed informally above.  A proof normalises when all detours can be eliminated. In intuitionistic logic, normalisation (i.e., that all proofs normalise) can be proven.\footnote{This is a central and standard result in proof theory. It was originally proved by \textcite{Prawitz1965-ao}. For an introductory presentation, see Chapter 6 of \textcite{Troelstra2000-ey} and Chapter 8 of \textcite{Negri2008-gp}} Prawitz interpreted normalisation as proving Gentzen correct that elimination rules use the main connective only in the sense afforded by the introduction rules, and this observation served as a foundation for his definition of proof-theoretic validity.    

We will present the formal definition of proof-theoretic validity in Section~\ref{sec:ptv}. For now, it is sufficient to note that an inference is proof-theoretically valid if it is one of the introduction rules for intuitionistic logic, or its occurrences in proofs are unnecessary. The second condition is illustrated by the reduction above which demonstrates how we can remove the rule of $\rightarrow$ elimination if we have a proof of the premises which ends in the rule of $\rightarrow$ introduction. That an inference is unnecessary is supposed to ensure that it follows from the introduction rules because in removing it, we demonstrate that, like the $\rightarrow$ elimination rule above, it only uses information provided by the introduction rules. \edit{It is worth noting that this conditions does not co-inside with the admissibility of a proof rule.\footnote{\edit{A proof rule is admissible in this context if given a closed proof of the premises, one can always find a closed proof of the conclusion.}} The elimination rules are not admissible in a system with only the introduction rules. The notion of admissibility is applied only to closed derivations, while proof-theoretic validity is defined on open and closed proofs.}

A logic is proof-theoretically valid if there are proof-theoretically valid proofs of all its theorems.  It is relatively easy to show that intuitionistic logic is proof-theoretically valid because this result follows from the proof-theoretic validity of intuitionistic elimination rules. Prawitz (\citeyear[246]{Prawitz1973-hc}) conjectured that intuitionistic logic was also the strongest proof-theoretically valid logic.  This conjecture remained open for almost thirty years before negative results started to appear (\cite{Sandqvist2009-dg}). \textcite{Piecha2019-sq} show that Prawitz's conjecture is false for all of the most important definitions of proof-theoretic validity. \edit{They show that superintuitionistic formulas, which we shall understand here to mean formulas that are classically consistent but not intuitionistically derivable, are proof-theoretically valid.}\footnote{\edit{It is important to note here that most notions of proof-theoretic validity aren't captured by a logic in the strict sense of a semi-decidable set of formulas closed under modus pones and substitution. However, we will use the term superintuitionistic to apply to any classically consistent collection of rules beyond intuitionistic logic.}} 

The superintuitionistic inference used to arrive at this result is Harrop's rule (\cite{Harrop1956-wk}):
$$\begin{prooftree}
\hypo{\neg A\rightarrow(B\vee C)}
\infer1{(\neg A\rightarrow B)\vee(\neg A\rightarrow C)}
\end{prooftree}$$
When added to intuitionistic logic, the result is Kreisel--Putnam logic. It turns out that this rule, or restrictions on it, are often proof-theoretically valid.\footnote{By `often' here I mean it holds in several prominent presentations of proof-theoretic validity (\cite{Piecha2015-ur, Piecha2019-sq}).  In \cite{Piecha2019-sq}, it is shown that for a very general class of systems that share properties with proof-theoretic validity either the system satisfies Harrop's rule, or Harrop's rule can be used to prove that the system must be superintuitionistic.}

It is worth highlighting that this result comes as a surprise since Prawitz claimed that the truth of his conjecture is obvious:
\begin{quote}
    It seems obvious that the elimination rules of Gentzen’s system are the elimination rules that correspond to his introduction rules. Or, again to put it more carefully: although there are of course weaker elimination rules and even elimination rules that are deductively equivalent with the ones formulated by Gentzen, there are no stronger rules that can be formulated in the language of predicate logic and are justifiable in terms of the introduction rules.\hfill\hbox{(\cite[270]{Prawitz2014-kk})}
\end{quote}
Prawitz's claim was supported by Dummett, who thought it `exceedingly plausible':
\begin{quote}
    It is exceedingly plausible that, on a verificationist meaning-theory, the correct logic will be intuitionistic;  \hfill\hbox{(\cite[270]{Dummett1991-mg})}
\end{quote}

The proof-theoretic validity of Harrop's rule is not just unexpected though: it threatens to undermine \edit{the motivation for} proof-theoretic validity.  For the moment, let us think of the introduction and elimination rules for a connective as \emph{harmonious} in case introduction of the connective followed by its elimination does not allow anything new to be proven. Then Harrop's rule is not harmonious \edit{because it isn't an elimination rule}.  If this is correct, we ought to conclude that proof-theoretic validity failed to offer a formalisation of Gentzen's idea.

I believe a defence can be mustered by carefully examining the causes of the non-intuitionistic inferences.  In sections \ref{sec:atomicformulas} and \ref{sec:happy} of this paper, I argue that the superintuitionistic inferences present in one of the most prominent notions of proof-theoretic validity stem from the treatment of atomic formulas and not from the treatment of connectives.  This means that, while proof-theoretic validity makes superintuitionistic rules valid, they do not function as elimination rules for the connectives.  Instead, it is more apt to think of them as being rules for atomic formulas.  We did not expect any rules to follow from the treatment of atomic formulas because it was thought that their treatment was entirely general.  It turns out, however, that it is not sufficiently general. Finally, if this analysis is correct then I propose a revised version of Prawitz's conjecture, which builds on a natural parallel between proof-theoretic validity and model-theoretic validity.  Intuitionistic logic is the strongest proof-theoretically valid logic if the definition of proof-theoretic validity is modified in such a way as to generalise the treatment of atomic formulas. \edit{Unlike Prawitz's original proposal, this conjecture proves to be correct ([citation removed to preserve blind review]).}

The paper is structured as follows. Section~\ref{sec:ptv} sets out the formal definition of proof-theoretic validity we will use. In Section~\ref{sec:atomicformulas}, we explore in detail a surprising complication in the definition of proof-theoretic validity that is caused by atomic formulas.  We discuss there how the inference rules used to define atomic formulas are equivalent to disjunction-free formulas.   In Section~\ref{sec:happy}, I then finally outline my argument to the effect that superintuitionistic inferences indeed follow from the treatment of atomic formulas.
\section{Proof-Theoretic Validity Defined}\label{sec:ptv}
As discussed in the introduction, proof-theoretic validity is an attempt to formally explicate what Gentzen meant when he claimed that elimination rules are consequences of the introduction rules. Central to proof-theoretic validity is normalisation and the reductions used in this process, which remove intuitionistic elimination rules. So, the validity of arguments is based in a sense on reducibility to the introduction rules.

While proof-theoretic validity is supposed to be a general notion, in the following we shall focus on the connectives `and', `or', and `if\dots then' with negation treated as $A\rightarrow\bot$.\footnote{Falsum is dealt with by always allowing the inference $\begin{prooftree}
    \hypo{\bot}\infer1{p}
\end{prooftree}$ for any atomic $p$.  We gloss over this as it will not impact the paper.}  Missing from this list are thus the quantifiers `all’ and `some’. If proof-theoretic validity without quantifiers leads to unharmonious validities, then adding quantifiers will not help but the addition of quantifiers would make the task of assessing the situation substantially more difficult. Therefore, in the following we limit ourselves to the propositional fragment.

To arrive at a definition of proof-theoretic validity, we will need to introduce the term \emph{argument} to refer to a potential proof. By argument, we mean something that has the right form to be a proof: a single conclusion that follows from assumptions, some of which may be discharged.\footnote{Formally, an argument is a tree labelled by formulas and a discharge relation.}  Proof-theoretic validity is a property of arguments, and arguments that satisfy the condition of proof-theoretically validity are proofs. Prawitz’s definition, as presented here, has four cases applying to arguments of different forms. 

To understand the four cases, we need to define what it means for an argument to be open or closed. An argument is closed if every top line of the argument is either an axiom or an assumption that has been discharged. An argument is open if it has an undischarged assumption. This gives us the cases: closed arguments that end in introduction rules, closed arguments that end in rules other than introduction rules, closed arguments that end in atomic formulas, and open arguments.

The simplest case in which an argument is valid is when it consists entirely of introduction rules. The introduction rules are supposed to define their main connectives, and an argument is valid when it follows from the introduction rules. There can thus be no problems with the use of introduction rules. When an introduction rule is used in an argument that contains other rules, then the application of the introduction rule is valid, while validity of the whole argument depends on validity of the rest of the argument. In the definition, this is accounted for by the following inductive case:
\begin{nmd}[Closed introduction case]
If $\mathcal{D}$ is a closed argument ending in an introduction rule, then it is valid if arguments for the premises of the introduction rule are valid.
\end{nmd}

When a closed argument ends in a rule that is not an introduction rule, we want it to be valid when it can be shown to follow, in some sense, from the introduction rules. Prawitz’s insight was that if an argument contains only intuitionistic introduction or elimination rules, then any closed argument ending in an elimination rule can be transformed via normalisation into a closed argument that ends in an introduction rule. In this way, one can step-by-step remove the elimination rules. This gives an initial condition for closed non-introduction arguments:
\begin{nmd}[Preliminary closed non-introduction case] 
If $\mathcal{D}$ is a closed argument that does not end in an introduction rule, then it is valid if \emph{it can be transformed by reductions used in the proof of normalisation into} a closed valid argument with the same conclusion that does end in an introduction rule.
\end{nmd}
\noindent Prawitz defends this definition by claiming that the argument with a detour and the argument that results from removal of the detour are syntactically distinct but semantically identical arguments (\cite[257]{Prawitz1971-ef};  \cite[234]{Prawitz1973-hc}). If two arguments are the same and one is valid (because of the use of introduction rules), then the other must also be valid.  

This method and justification can only work when one considers proofs in intuitionistic logic, not arbitrary arguments. But we must consider arbitrary arguments, otherwise, we have stacked the deck in favour of intuitionistic logic. To resolve this problem, we need to generalise the concept of normalisation to arbitrary arguments. Prawitz does this by arguing that there must be a computable transformation of an argument to a valid argument with the same conclusion that ends in an introduction rule  (\cite{Prawitz1973-hc}). \textcite[552-3]{Schroeder-Heister2006-rh} does away with the requirement that the transformation be computable. No reason is given, in either case, to think that the transformed argument is semantically identical to the original.  

As presented here, the inductive case does not deal with transformations on arguments, but given a sufficiently broad understanding of transformations, it can be shown that the definition is equivalent to one that does. Not including reference to transformations is therefore merely a technical convenience. Thus, what we end up with is the following condition:
\begin{nmd}[Closed non-introduction case] 
If $\mathcal{D}$ is a closed argument that does not end in an introduction rule, it is valid if \emph{there is} a closed valid argument with the same conclusion that does end in an introduction rule.
\end{nmd}
\noindent With these two conditions on closed arguments in place, we can move on to arguments with open assumptions.

An open assumption is replaced with a valid argument (or proof). This demonstrates a key commitment embedded in the idea that a connective can be defined by a proof rule. When a sentence with a logical connective is asserted, it is taken to be a claim that there is a valid argument (or proof) of the sentence. An assumption is then an assertion with its assertoric force cancelled. Reasoning that follows from an assumption is hypothetical on the statement’s being asserted. As such, the reasoning that follows is valid if it would be valid were there a valid argument (or proof) of the assumptions. It is this fundamental idea motivating the treatment of open assumptions: for an argument with open assumptions to be valid, the reasoning must be valid if we had valid arguments (or proofs) of the assumptions. It does not suffice, however, to check just one valid argument (or proof): any closed valid arguments (or proofs) substituted for assumptions must result in a valid argument. This gives us the following preliminary condition:
\begin{nmd}[Preliminary open case] If $\mathcal{D}$ is an open argument of $ A $ with open assumptions $ A _0,\dots, A _n$, then it is valid if for all closed valid arguments $\mathcal{D}_0,\dots,\mathcal{D}_n$ of $ A _0,\dots, A _n$, the following argument is valid:
        \begin{center}
            \begin{prooftree}
            \hypo{\mathcal{D}_0}
            \infer[no rule]1{ A _0}
            \hypo{ \phantom{A}$\dots$\phantom{A}}
            \infer[no rule]1{ \phantom{A}$\dots$\phantom{A}}
            \hypo{\mathcal{D}_n}
            \infer[no rule]1{ A _n}
            \infer[no rule]3{\mathcal{D}}
            \infer[no rule]1{ A }
        \end{prooftree}
        \end{center}
\end{nmd}
However, our treatment of arguments with open assumptions is currently incomplete.  The reason for this is that as it stands, we have not explained how to deal with arguments such as the following, where $p$, $q$, and $r$ are atomic formulas:\footnote{We shall use lowercase letters for atomic formulas and uppercase letters for arbitrary formulas.  It is important in this setting to distinguish between the two.}
\begin{ex}\label{ex:atomicmp}
\begin{center}
    \begin{prooftree}
    \hypo{p\rightarrow q}
    \hypo{p}
    \infer2{q}
    \end{prooftree}
\end{center}
\end{ex}
\noindent This is because there are no closed valid arguments (or proofs) of the atomic formula $p$ as things have been set up so far.  None of the intuitionistic introduction rules end in an atomic formula, and so no valid argument (or proof) of $p$ thus ends in one of the intuitionistic introduction rules.  Note that $p$ as an arbitrary propositional variable stands here for any atomic sentence.  In model-theoretic semantics, arbitrary propositions are dealt with by considering the different meanings they could have, whereby meanings are spelt out in terms of truth values.  Here, we take terms to be defined by their proof rules; in analogy with valuations in model-theoretic semantics, we are thus going to consider different assignments of proof rules to the atomic sentences.  Let us relativise validity to a set $S$ of rules for atomic formulas. We will say a lot more about what counts as an atomic rule in the following section.\footnote{For now, one can think of this as containing atomic axioms, such as $p$ or $q$, and inferences from atomics to atomics, such as $\frac{p}{q}$ or $\frac{r \ s \ q}{p}$.}  This can be understood as something like a theory that constrains the use--and hence the interpretation--of the propositional letters.  Then we require what is valid to be relativised to a particular assignment of meaning to the atomic sentences. 

Before revising our definition for open arguments, we can now give the condition for atomic formulas:
\begin{nmd}[Atomic case] If $\mathcal{D}$ is a closed argument ending in an atomic formula, then it is $S$-valid if there is a proof of $p$ containing only rules in $S$.
\end{nmd} 
\noindent This counts as the base case of the definition because when we have an atomic formula, we can restrict our attention to arguments that do not contain logical connectives. 

Initially, we took this detour because we did not know what to do when we had an open assumption of an atomic formula. We see now that we will need to substitute the assumption for a valid argument (or proof) of the atomic formula, and we know what such a valid argument (or proof) looks like by the condition spelt out above. Because it is relative to a set of atomic rules, we will need to make all our other conditions relative to a set of atomic rules by replacing validity with  $S$-validity in the definitions.  But for open assumptions, we are going to make one other change.  

Consider Example~\ref{ex:atomicmp} again: what happens with $S$-validity when $S$ is a set of atomic rules from which $p$ is not provable?  In this case, our preliminary open case is trivially satisfied because there are no valid arguments for $p$.  But this is not what we want to happen; our reasoning is hypothetical on there being a proof for $p$.\footnote{There is a variant of this notion where things are still relative to a set of atomic rules $S$ but no extensions are considered.  This at one point was Prawitz's preferred notion (\cite{Prawitz2014-kk}).  Both this and the displayed definition are superintuitionistic (\cite[thm 2.5]{Piecha2019-sq}).} To get around this, we will consider not just our initial set of atomic rules $S$ but also any extensions of it, and we can extend $S$ so that $p$ is provable.  Hopefully, this will mean that for every atomic $p$, we can find an extension of $S$ in which there is a valid argument for $p$. This gives us the condition for arguments with open assumptions:

\begin{nmd}[Open case] If $\mathcal{D}$ is an open argument of $ A $ with open assumptions $ A _0,\dots, A _n$, it is $S$-valid if for all $S'$ which are acceptable extensions of $S$ and all closed $S'$-valid arguments $\mathcal{D}_0,\dots,\mathcal{D}_n$ of $ A _0,\dots, A _n$, the following argument is $S'$-valid:
        \begin{center}
            \begin{prooftree}
            \hypo{\mathcal{D}_0}
            \infer[no rule]1{ A _0}
            \hypo{ \phantom{A}$\dots$\phantom{A}}
            \infer[no rule]1{ \phantom{A}$\dots$\phantom{A}}
            \hypo{\mathcal{D}_n}
            \infer[no rule]1{ A _n}
            \infer[no rule]3{\mathcal{D}}
            \infer[no rule]1{ A }
        \end{prooftree}
        \end{center}
\end{nmd}
This gives us the following definition:

\begin{defn}(\cite[236]{Prawitz1973-hc}; \cite[543-4]{Schroeder-Heister2006-rh})\label{def:PTvalid} An argument $\mathcal{D}$ is \emph{$S$-valid} for a set of rules $S$ describing the behaviour of the atomic formulas if one of the following conditions holds:
    \begin{enumerate}[wide, labelwidth=!, labelindent=0pt]
        \item[{\bf Atomic case}] If $\mathcal{D}$ is a closed argument ending in an atomic formula, then it is $S$-valid if there is a proof of $p$ containing only rules in $S$.
        \item[\bf Closed introduction case] If $\mathcal{D}$ is a closed argument ending in an introduction rule, then it is $S$-valid if arguments for the premises of the introduction rule are $S$-valid.
        \item[\bf Closed non-introductory case] If $\mathcal{D}$ is a closed argument which does not end in an introduction rule, then it is $S$-valid if there is a $S$-valid argument with the same conclusion that does end in an introduction rule.  
        \item[\bf Open case] If $\mathcal{D}$ is an open argument of $ A $ with open assumptions $ A _0,\dots, A _n$, it is $S$-valid if for all $S'$ which are acceptable extensions of $S$ and all closed $S'$-valid arguments $\mathcal{D}_0,\dots,\mathcal{D}_n$ of $ A _0,\dots, A _n$, the following argument is $S'$-valid:
        \begin{center}
            \begin{prooftree}
            \hypo{\mathcal{D}_0}
            \infer[no rule]1{ A _0}
            \hypo{ \phantom{A}$\dots$\phantom{A}}
            \infer[no rule]1{ \phantom{A}$\dots$\phantom{A}}
            \hypo{\mathcal{D}_n}
            \infer[no rule]1{ A _n}
            \infer[no rule]3{\mathcal{D}}
            \infer[no rule]1{ A }
        \end{prooftree}
        \end{center}
    \end{enumerate}
\end{defn}

A rule of inference is then \emph{proof-theoretically valid} if it is valid on all acceptable sets of atomic rules. That is $S$-valid for all $S$. We can write this as a consequence relation as follows $\vDash_{PTV}\varphi\Leftrightarrow $ for all acceptable set of atomic rules $S$, there is a $S$-valid proof of $\varphi$.  But this definition is imprecise because we have not said what an acceptable set of atomic rules is.  It turns out that fixing this issue leads to a zoo of distinct versions of proof-theoretic validity, \edit{all of which contain superintuitionistic validates.}
\section{Atomic Rules and Proof-Theoretic Systems}\label{sec:atomicformulas}
In the previous section, I laid out the definition of proof-theoretic validity (or at least one version of it). Here we will discuss how the set of atomic rules used in the definition affects what is proof-theoretically valid, as has been shown by a growing body of formal work (\cite{Sandqvist2009-dg}, \cite{Piecha2015-ur}, \cite{Goldfarb2016-tk}, \cite{Piecha2019-sq}, \cite{Stafford}).

To discuss this in detail, we first need to discuss what an atomic rule is. The easiest rules to understand are those that look like axioms. Axioms are immediate inferences that have a conclusion but no premises. Axioms can be written down as follows: $\begin{prooftree}\hypo{}\infer1{p}\end{prooftree}$.  We can, however, also have atomic rules that look like the proof rules for connectives but contain only atomic formulas.  The following example contrasts an atomic rule with $\wedge$ introduction.
\begin{ex}\label{ex:level1rule}
$$\begin{prooftree}\hypo{ A }\hypo{ B }\infer2{ A \wedge B }\end{prooftree}\qquad\begin{prooftree}\hypo{p}\hypo{q}\infer2{r}\end{prooftree}$$
\end{ex}
\noindent Similar but slightly more complicated would be the case of a rule that discharges assumptions such as $\rightarrow$ introduction or $\vee$ elimination.
\begin{ex}\label{ex:level2rule}
$$\text{\bf a. \ }\begin{prooftree}\hypo{[ A ]}\ellipsis{}{ B }\infer1{ A \rightarrow B }\end{prooftree}\qquad\begin{prooftree}\hypo{[p]}\ellipsis{}{q}\infer1{r}\end{prooftree}\qquad\text{\bf b. \ }\begin{prooftree}\hypo{ A \vee B }\hypo{[ A ]}\ellipsis{}{ C }\hypo{[ B ]}\ellipsis{}{ C }\infer3{ C }\end{prooftree}\qquad\begin{prooftree}\hypo{p}\hypo{[q]}\ellipsis{}{s}\hypo{[r]}\ellipsis{}{s}\infer3{s}\end{prooftree}$$
\end{ex}

\edit{It is worth noting at this point that we need to withhold the common instinct to read sentences with atomic formulas such as  $(p\wedge q)\rightarrow r$ as just another way of writing $( A \wedge B )\rightarrow C $.  When an atomic letter such as $p$ is used, it should be understood as an atomic sentence. When a schematic letter $A$ is used, it should be understood as an arbitrary sentence and one can substitute in anything for $A$ unless otherwise specified.}  

Recall that in the definition of proof-theoretic validity, one of the cases involved \emph{acceptable extensions} of a set of rules $S$.  We might take the sets of atomic rules we consider in the definition of proof-theoretic validity to be any collection of rules of the three types (axioms, premise-to-conclusion, and with assumptions) we have just discussed. In that case, an acceptable extension would be any set that adds rules of these types. We could, however, imagine that one might reject atomic rules which discharge hypotheses as acceptable atomic rules: in such case, one would arrive at a different definition of what an acceptable extension is.

Still, either definition of an acceptable extension would yield a set of validates so strange they do not deserve the name logic. To see this, we need to note that there is a relationship between these atomic rules and formulas. First, note that we can treat the inference from premises to the conclusion as an implication,  $\rightarrow$, and we simply conjoin the premises.  So $\begin{prooftree}\hypo{p}\hypo{q}\infer2{r}\end{prooftree}$ can be transformed into $(p\wedge q)\rightarrow r$.  The following two proofs show that in arguments, the formula and the rule are interchangeable. The parts that correspond to the rule are highlighted in bold to make the proofs easier to read.
$$\begin{prooftree}\hypo{[p\wedge q]}\infer1{\bf p}\hypo{[p\wedge q]}\infer1{\bf q}\infer2{\bf r}\infer1{(p\wedge q)\rightarrow r}\end{prooftree}\qquad\begin{prooftree}\hypo{\bf p}\hypo{\bf q}\infer2{p\wedge q}\hypo{(p\wedge q)\rightarrow r}\infer2{\bf r}\end{prooftree}$$
The first proof shows that the formula can be proven using the rule and the second proof shows that the rule can be, in an argument, replaced by the formula. It can be proven that any rule from atomic premises to atomic conclusion is equivalent to a formula in just this way (\cite{Piecha2015-ur}).

This can also be done when rules discharge assumptions. We just need to treat the proof from the assumption to its premise as an implication.  Example~\ref{ex:level2rule}.a is thus interchangeable with $(p\rightarrow q)\rightarrow r$ and Example~\ref{ex:level2rule}.b is interchangeable with $(p\wedge(q\rightarrow s)\wedge(r\rightarrow s))\rightarrow s$.  Using Example~\ref{ex:level2rule}.b, we can prove this interchangeability as follows:\footnote{Note that I am using a slight generalisation of $\wedge$ introduction and elimination.  This is only to make the proofs more perspicuous and can be easily removed at no cost.}

{$$\begin{prooftree}\hypo{[p\wedge(q\rightarrow s)\wedge(r\rightarrow s)]}\infer1{{\bf p}}\hypo{\bf[q]}\hypo{[p\wedge(q\rightarrow s)\wedge(r\rightarrow s)]}\infer1{q\rightarrow s}\infer2{\bf s}\hypo{\bf [r]}\hypo{[p\wedge(q\rightarrow s)\wedge(r\rightarrow s)]}\infer1{r\rightarrow s}\infer2{\bf s}\infer3{\bf s}\infer1{(p\wedge(q\rightarrow s)\wedge(r\rightarrow s))\rightarrow s}\end{prooftree}$$}
\vspace{.2cm}

$$\begin{prooftree}\hypo{\bf p}\hypo{\bf [q]}\ellipsis{}{\bf s}\infer1{q\rightarrow s}\hypo{\bf [r]}\ellipsis{}{\bf s}\infer1{r\rightarrow s}\infer3{p\wedge(q\rightarrow s)\wedge(r\rightarrow s)}\hypo{(p\wedge(q\rightarrow s)\wedge(r\rightarrow s))\rightarrow s}\infer2{\bf s}\end{prooftree}$$

If we use these three different types of rules, some--but not all--formulas are interchangeable with rules.\footnote{Specifically, formulas built up from $\wedge$ and $\rightarrow$ where no branch of the syntax tree contains more than two $\rightarrow$s are equivalent to the rules.} Because of this, if we take acceptable sets of atomic rules in the definition of proof-theoretic validity to include only atomic rules of these types, an odd thing will happen.  An intuitionistically invalid rule called generalised Harrop's rule will be proof-theoretically valid,

\begin{equation}\tag{Generalised Harrop's rule}
    \begin{prooftree}
    \hypo{ A \rightarrow( B \vee C )}
    \infer1{( A \rightarrow B )\vee( A \rightarrow C )}
    \end{prooftree}
\end{equation}

\noindent but only when $A$ is exchanged for formulas interchangeable with rules (\cite{Piecha2015-ur}). Note that Harrop's rule is generalised Harrop's rule with $A$ restricted to formulas beginning with a negation.

Let us look at why \edit{cases of} generalised Harrop's rule is valid using an example.  We will make use of the equivalence with formulas and all that is required for the proof is that $S\cup\set{\begin{prooftree}
    \hypo{}\infer1{p}
\end{prooftree}}$ is an acceptable set of atomic rules if $S$ is.  
\begin{ex}\label{ex:harrop}$$
\begin{prooftree}
    \hypo{p\rightarrow(q\vee r)}
    \infer1{(p\rightarrow q)\vee(p\rightarrow r)}
\end{prooftree}
$$\end{ex}
It turns out that asking whether the above inference is $S$-valid is--by the open case of Definition~\ref{def:PTvalid}--the same thing as asking whether for all $S'$, extending $S$, when $p\rightarrow(q\vee r)$ is $S'$ valid so is $(p\rightarrow q)\vee(p\rightarrow r)$. That in turn, via the above-mentioned equivalence between rules and formulas, can be answered by looking whether when $q\vee r$ is $S'\cup\set{\Bar{p}}$-valid, then $(p\rightarrow q)\vee(p\rightarrow r)$ is $S'$-valid.   But if $q\vee r$ is $S'\cup \set{\Bar{p}}$-valid, then by a combination of the closed introduction case and closed non-introduction case, there is a $S'\cup\set{\Bar{p}}$-valid argument for $q\vee r$ that ends in $\vee$-introduction. There is thus a $S'\cup \set{\Bar{p}}$-valid argument for either $q$ or $r$.   Assume $q$ is $S'\cup \set{\Bar{p}}$-valid, then again--by the equivalence between rules and formulas--it follows that $p\rightarrow q$ is $S'$-valid, and by the closed introduction case, so is $(p\rightarrow q)\vee(p\rightarrow r)$.

Let us take a moment to restate the situation.  We now have the first precise definition of proof-theoretic validity provided by filling in `acceptable set of atomic rules' with `set of atomic rules containing axioms, premises-to-conclusion and with assumptions'.  When we do this, instances of Harrop's rule are proof-theoretically valid.  

Having the generalised Harrop’s rule hold for a rather circumscribed and arbitrary collection of sentences is not ideal but this issue can be partially resolved by expanding the equivalence between formulas and rules.   We saw that $(p\rightarrow q)\rightarrow r$ corresponds to the rule in Example~\ref{ex:level2rule}.a. What would then correspond to a formula such as $((p\rightarrow q)\rightarrow r)\rightarrow s$? \textcite{Schroeder-Heister1984-mx} proposed a generalisation of inference rules that allows us to find a corresponding rule for this formula, but it requires an expansion of what can be discharged.  Note that in the formula corresponding to the discharge of an assumption, we had implication in two places: the first standing in for the proof from the assumption to the premise and the second standing for the inference from premise to conclusion.  To add an extra implication, we need to add another of these.  We will do this by assuming not an atomic formula, but a rule.  The resulting rule is: 
\begin{ex}\label{ex:lv3}
$$\begin{prooftree}
        \hypo{p}
        \infer1{q}
        \delims{\left[}{\right]}
        \infer[no rule]1{\vdots}
        \infer[no rule]1{r}
        \infer1{s}
    \end{prooftree}$$
\end{ex}
\noindent where the square brackets represent discharge of the assumed rule.\footnote{\edit{Note that, while all our examples have the rule discharged at the top of the proof, the rule discharged can appear anywhere in the proof tree.}} We can see this rule is equivalent to $((p\rightarrow q)\rightarrow r)\rightarrow s$ via the following two proofs:
$$\begin{prooftree}\hypo{[(p\rightarrow q)\rightarrow r]^3}\hypo{\bf [p]^1}\infer1{\bf q}
        \delims{\left[}{\right]^2}\infer1[1]{p\rightarrow q}\infer2{\bf r}
\infer1[2]{\bf s}
\infer1[3]{((p\rightarrow q)\rightarrow r)\rightarrow s}\end{prooftree}\qquad\qquad \qquad\begin{prooftree}\hypo{}
        \ellipsis{}{\bf p}\hypo{[p\rightarrow q]^1}\infer2{\bf q}\ellipsis{}{\bf r}\infer1[1]{(p\rightarrow q)\rightarrow r}\hypo{((p\rightarrow q)\rightarrow r)\rightarrow s}\infer2{\bf s}\end{prooftree}$$

Let us look quickly at how we would use a rule like this in a proof.  First, let us construct a proof using the rule in Example~\ref{ex:lv3} along with the axiom $\begin{prooftree}
\hypo{}\infer1{p}
\end{prooftree}$ and the one step inference $\begin{prooftree}
\hypo{q}\infer1{r}
\end{prooftree}$
\begin{ex}
$$\begin{prooftree}
    \hypo{\bar{p}}
    \infer1{q}
    \delims{\left[}{\right]^1}
    \infer1{r}
    \infer1[$_1$]{s}
\end{prooftree}$$
\end{ex}
\noindent This can be generalised even further.  We can discharge rules that discharge assumptions themselves and in doing so acquire more equivalences between formulas and rules.  As a second example, take the following proof of $t$ with a rule that lets you discharge the rule in Example~\ref{ex:lv3} when going from $s$ to $t$. Note that the first four lines are identical to the proof above.
\begin{ex}
$$
\begin{prooftree}
    \hypo{\bar{p}}
    \infer1{q}
    \delims{\left[}{\right]^1}
    \infer1{r}
    \infer1[$^1$]{s}
    \delims{\left[}{\right]^2}
    \infer1[$_2$]{t}
\end{prooftree}$$
\end{ex}

If we allow rules that discharge rules of any level of complexity, then every disjunction-free formula is equivalent to a rule (\cite{Piecha2015-ur}).  This means that we have generalised Harrop's formula when the antecedent $A$ is disjunction-free. \edit{Then the set of proof-theoretically valid formulas is the same as the set of consequences of generalised inquisitive logic, i.e., intuitionistic logic plus the generalised Harrop's rule restricted to disjunction-free formulas (\cite{Piecha2015-ur,Stafford}). This logic is far more uniform than the set of validities resulting from considering only atomic rules that do not discharge rules, but it is still not closed under uniform substitution because we do not have Harrop's rule for all disjunctive formulas.}\footnote{\edit{We use the term logic as that is what is used in the existing literature. However, there is a general view that logics should be closed under subsitution. This view can be found in early works on modern logic such as Bolzano:
\begin{quote}
    But suppose that there is just a single idea in it which can be arbitrarily varied without disturbing its truth or falsity, i.e. if all the propositions produced by substituting for this idea any other idea we pleased are either true altogether or false altogether, presupposing that they have a denotation. [\ldots] I permit myself, then, to call propositions of this kind, borrowing an expression from Kant, \emph{analytic}.\hfill\hbox{(\cite[\S{148} p. 192]{Bolzano1973-jq})}
\end{quote}
The works of Tarski:
\begin{quote}Consider any class $K$ of sentences and a sentence $X$ which follows from the sentences of this class. [\dots] 
Moreover, since we are concerned here with the concept of logical, i.e. \emph{formal}, consequence, and thus with a relation which is to be uniquely determined by the form of the sentences between which it folds, this relation cannot be influenced in any way by empirical knowledge, and in particular by knowledge of the objects to which the sentence $X$ or the sentences of the class $K$ refer.
The consequence relation cannot be affected by replacing the designations of the objects referred to in these sentences by the designations of any other objects.\hfill\hbox{(\cite[414-5]{Tarski1956-um})}\end{quote}
And more modern works for instance by Stalnaker:
\begin{quote}
    A substitution instance of an argument is an argument of the same form, and arguments are logically valid only if all arguments of the same form are logically valid. \hfill\hbox{(\cite[332]{Stalnaker1977-rw})}
\end{quote}
However, there are other notions of uniform substitution.  For example, Humberstone offers the following:
\begin{quote}
    [...] closure under uniform substitution of propositional variables (rather than arbitrary formulas) for propositional variables.\hfill\hbox{(\cite[188]{Humberstone2011-sd})}
\end{quote}}
}

We see above that we have choices of how to treat the atomic formulas.  We could consider proof-theoretic validity with only atomic rules that do not discharge anything or we could allow any atomic rules even those that discharge other rules.  It will be useful to have explicit notation for what sets of atomic rules are allowed. We shall call any set of sets of atomic rules a \emph{proof-theoretic system}.   If we allow any possible set of atomic rules lets use the shorthand the {\it complete system} and if we don't allow atomic rules to discharge assumptions we will write the {\it minimal system}.   We can have more esoteric systems as well, such as $\set{S\mid S\text{ is in the complete system and }\Bar{p}\notin S}$ or even $\set{\varnothing, \set{\Bar{p}},\set{\Bar{q}}}$.  Each of these will give us a different logic or theory which is proof-theoretically valid.  So, the question of whether a derivation $\mathcal{D}$ is proof-theoretically valid will be relative to a proof-theoretic system. We can modify the proof-theoretic consequence relation so that if $\mathfrak{S}$ is a proof-theoretic system then $\vDash_{PTV\mathfrak{S}}\varphi$ if and only if for all $S\in\mathfrak{S}$, there is an $S$-valid proof of $\varphi$ where all acceptable extensions of $S$ are in $\mathfrak{S}$.

\section{The Role of the Atomic Formulas Explained}\label{sec:happy}
In this section, it will be argued that superintuitionistic \edit{rules} that are proof-theoretically valid follow from the treatment of atomic formulas and not, as one might initially think, from the treatment of the connectives. As such, the superintuitionistic validates should not be viewed as undermining the argument – partially supported by proof-theoretic validity – that intuitionistic elimination rules are privileged. Finally, by drawing a comparison with model-theoretic semantics, a modified Prawitz's conjecture will be proposed \edit{that is provably correct}.

We will consider the validity of the one-step inferences    \scalebox{0.7}{$
\begin{prooftree}
    \hypo{p}
    \infer1{q\vee r}
\end{prooftree}
$}. 
Its validity is anticorrelated with the validity of Example~\ref{ex:harrop}. \edit{The discussion in section~\ref{sec:atomicformulas} of the validity Example~\ref{ex:harrop} can be extended to demonstrate the anticorrolation but we leave this to the reader.} We will show that \scalebox{0.7}{$
\begin{prooftree}
    \hypo{p}
    \infer1{q\vee r}
\end{prooftree}
$} is valid in the proof-theoretic system $\set{{\tiny\varnothing,\set{\begin{prooftree}\hypo{}\infer1{p}\end{prooftree},\begin{prooftree}\hypo{p}\infer1{q}\end{prooftree}},\set{\begin{prooftree}\hypo{}\infer1{p}\end{prooftree},\begin{prooftree}\hypo{p}\infer1{r}\end{prooftree}}}}$, which is in some sense a toy system, but will do nicely for the purposes of illustration.

It is sufficient to show that \scalebox{0.7}{$
\begin{prooftree}
    \hypo{p}
    \infer1{q\vee r}
\end{prooftree}
$}
\emph{is} $\varnothing$-valid.  This follows, by the condition for open arguments, if for every extension $S$ of $\varnothing$ and for every argument for $p$ that is $S$-valid, the result of appending the argument for $p$ to \scalebox{0.7}{$
\begin{prooftree}
    \hypo{p}
    \infer1{q\vee r}
\end{prooftree}
$} is $S$-valid.  As it happens in our small system, we only have one valid argument (or proof) of $p$, namely  
$
\begin{prooftree}
    \hypo{}
    \infer1{p}
\end{prooftree}
$.
This argument is both $\set{{\tiny\begin{prooftree}\hypo{}\infer1{p}\end{prooftree},\begin{prooftree}\hypo{p}\infer1{q}\end{prooftree}}}$ and $\set{{\tiny\begin{prooftree}\hypo{}\infer1{p}\end{prooftree},\begin{prooftree}\hypo{p}\infer1{r}\end{prooftree}}}$-valid by the condition for the atomic case.  So, it needs to be shown that
$$
\begin{prooftree}
    \hypo{}
    \infer1{p}
    \infer1{q\vee r}
\end{prooftree}
$$
is both $\set{{\tiny\begin{prooftree}\hypo{}\infer1{p}\end{prooftree},\begin{prooftree}\hypo{p}\infer1{q}\end{prooftree}}}$ and $\set{{\tiny\begin{prooftree}\hypo{}\infer1{p}\end{prooftree},\begin{prooftree}\hypo{p}\infer1{r}\end{prooftree}}}$-valid. But this follows in the first case because 
$$
\begin{prooftree}
    \hypo{}
    \infer1{p}
    \infer1{q}
    \infer1{q\vee r}
\end{prooftree}
$$
is $\set{{\tiny\begin{prooftree}\hypo{}\infer1{p}\end{prooftree},\begin{prooftree}\hypo{p}\infer1{q}\end{prooftree}}}$-valid, by the condition for the closed introduction case, and in the second case we can simply replace $q$ with $r$ in the above proof and it will likewise be $\set{{\tiny\begin{prooftree}\hypo{}\infer1{p}\end{prooftree},\begin{prooftree}\hypo{p}\infer1{r}\end{prooftree}}}$-valid.


A small change in the proof-theroetic system will lead to \scalebox{0.7}{$
\begin{prooftree}
    \hypo{p}
    \infer1{q\vee r}
\end{prooftree}
$} being invalid. Take the previous system $\set{{\tiny\varnothing,\set{\begin{prooftree}\hypo{}\infer1{p}\end{prooftree},\begin{prooftree}\hypo{p}\infer1{q}\end{prooftree}},\set{\begin{prooftree}\hypo{}\infer1{p}\end{prooftree},\begin{prooftree}\hypo{p}\infer1{r}\end{prooftree}}}}$ and add $\set{\begin{prooftree}\hypo{}\infer1{p}\end{prooftree}}$, the obtaining the system $\set{{\tiny\varnothing,\set{\begin{prooftree}\hypo{}\infer1{p}\end{prooftree}},\set{\begin{prooftree}\hypo{}\infer1{p}\end{prooftree},\begin{prooftree}\hypo{p}\infer1{q}\end{prooftree}},\set{\begin{prooftree}\hypo{}\infer1{p}\end{prooftree},\begin{prooftree}\hypo{p}\infer1{r}\end{prooftree}}}}$. To test the $\varnothing$-validity of 
\scalebox{0.7}{$
\begin{prooftree}
    \hypo{p}
    \infer1{q\vee r}
\end{prooftree}
$} in this system, we need to check whether every extension of $\varnothing$ that has a valid argument (or proof) of $p$ has a valid argument (or proof) of $q\vee r$. Only one thing has changed from our earlier discussion, namely that $
\begin{prooftree}
    \hypo{}
    \infer1{p}
\end{prooftree}
$ is $\set{\begin{prooftree}\hypo{}\infer1{p}\end{prooftree}}$-valid. But we cannot give $\set{\begin{prooftree}\hypo{}\infer1{p}\end{prooftree}}$-valid argument (or proof) of $q$ or $r$ by the condition for atomic conclusions. It follows from this that \scalebox{0.7}{$
\begin{prooftree}
    \hypo{p}
    \infer1{q\vee r}
\end{prooftree}
$} is not $\varnothing$-valid in this proof-theoretic system. 

This is a toy example, but it mimics what happens when working with larger proof-theoretic systems. When a system lacks a rule, some one-step inferences from a formula to a disjunction will be valid, e.g.,\@ from $A$ to $B\vee C$, but by adding an atomic rule that is equivalent to $A$, this one-step inference becomes invalid.  This is because there is now an additional set of atomic rules that has a valid argument (or proof) of $A$ (because there is a new rule that is equivalent to it) but not of either $B$ or $C$. As the validity of generalised Harrop's rule realise on the invalidity of the one step inference from antecedent to disjunctive consequence, I take this to show that the validity of generalised Harrop's rule follows, not from the treatment of connectives, but rather from the correspondence between atomic rules and formulas, and ultimately from the treatment of atomic formulas. The treatment of the atomic formulas in the complete system has caused disjunction-free formulas to have no relationship to one another, even only potentially.  As a result, generalised Harrop's rule follows for those formulas.

One might think that to solve this issue, we should just remove the equivalence between formulas and rules to get a system in which the treatment of atomic formulas does not have this impact.  That, however, will not work for two reasons.  The first is that even with the simplest of rules for atomic formulas, the axioms, there is still an equivalence between the axiom $\begin{prooftree}\hypo{}\infer1{p}\end{prooftree}$ and the formula $p$.  And a system with no atomic rules would behave very strangely indeed.\footnote{One could change the treatment of open formulas as well and there are versions of proof-theoretic validity that behave like this.  They have not been more successful in capturing intuitionistic logic, however.}  Without atomic rules, $p\rightarrow q$ would have a closed valid argument (or proof) because there are no closed valid arguments (or proofs) of $p$, and any closed valid argument (or proof) of $p$ is thus also one of $q$.  The second is that it is not just this treatment of the atomic formulas that causes the problem: it is any treatment.  It has been shown that for all common notions, \edit{the resulting set of validates will be superintuitionistic} (\cite{Piecha2019-sq}).  For example, if we take the minimal system, i.e., one where no atomic rules can discharge assumptions, it turns out that $\neg\neg p\rightarrow p$ is valid (\cite{Sandqvist2009-dg}).\footnote{Remember that these \edit{systems} are not closed under substitutions so this does not imply that $\neg\neg A\rightarrow A$ is valid which would mean the \edit{system} was classical \edit{logic}.}

What is needed instead is a more general, not more restrictive, treatment of the atomic formulas.  Recalling our two toy examples, we need a notion of proof-theoretic validity that allows for both the possibilities that they represent.  That is, $p$ may imply $q\vee r$ or it may not.  We cannot achieve this by creating a larger proof-theoretic system.  After all, the union of $\set{{\tiny \varnothing,\set{\begin{prooftree}\hypo{}\infer1{p}\end{prooftree},\begin{prooftree}\hypo{p}\infer1{q}\end{prooftree}},\set{\begin{prooftree}\hypo{}\infer1{p}\end{prooftree},\begin{prooftree}\hypo{p}\infer1{r}\end{prooftree}}}}$ and $\set{{\tiny\varnothing,\set{\begin{prooftree}\hypo{}\infer1{p}\end{prooftree}},\set{\begin{prooftree}\hypo{}\infer1{p}\end{prooftree},\begin{prooftree}\hypo{p}\infer1{q}\end{prooftree}},\set{\begin{prooftree}\hypo{}\infer1{p}\end{prooftree},\begin{prooftree}\hypo{p}\infer1{r}\end{prooftree}}}}$ is just $\set{{\tiny\varnothing,\set{\begin{prooftree}\hypo{}\infer1{p}\end{prooftree}},\set{\begin{prooftree}\hypo{}\infer1{p}\end{prooftree},\begin{prooftree}\hypo{p}\infer1{q}\end{prooftree}},\set{\begin{prooftree}\hypo{}\infer1{p}\end{prooftree},\begin{prooftree}\hypo{p}\infer1{r}\end{prooftree}}}}$  again. 

To see how we can do this, I would like to draw a parallel with Tarskian model-theoretic semantics.  In Tarskian model-theoretic semantics for classical logic, classical validities are those which hold in all models.  But there is no one model in which all and only the classical validities hold.  In fact, despite no atomic formula $p$ or its negation $\neg p$ being a classical validity, \edit{every classical model models} either $p$ or $\neg p$ for all atomic $p$.  Does this mean that these superclassical validities do not follow from the treatment of atomic formulas?  Clearly not: it is to be expected that every model will have many superclassical validities.  I would suggest that proof-theoretic semantics can use a similar approach.  While every proof-theoretic system has some superintuitionistic validities, to ask what is proof-theoretically valid in general we must look at \emph{all} proof-theoretic systems.  

\edit{We can summerise this line of argument as follows: When a choice is made about what counts as an extension of a set of atomic rules, that choice necessarily leads to some superintuitionistic validities. But it turns out that for every superintuitionistic validity there is a choice of what counts as extensions of the atomic rules which is a counterexample to it.  We can express this claim formally as follows: if $\mathbb{S}$ is the set of all atomic rules, then for any $\mathfrak{S}\subseteq\mathcal{P}(\mathbb{S})$ there is a formula $\varphi$ not provable in intuitionistic logic such that $\mathfrak{S}\vDash\varphi$ (\cite{Piecha2019-sq}). But also for any superintuitionistic formula $\varphi$ there is a $\mathfrak{S}\subseteq\mathcal{P}(\mathbb{S})$ such that $\mathfrak{S}\nvDash\varphi$ ([citation removed to preserve blind review]). 

Prawitz’s conjecture as initially posed was about a particular choice of what counts as an extension of the set of atomic rules. And what \textcite{Piecha2019-sq} showed is that there is no way of making this choice such that Prawitz’s conjecture holds. While this conclusively refutes the technical question that had been unresolved since the 70s, here I propose a modification of the conjecture. The discussion above has aimed to demonstrate that there is no choice of what counts as an extension of a set of atomic rules that does not encode information about the atomic rules into the notion of validity. As such, if we wish to have our logic be neutral to the non-logical content we should instead define what is proof-theoretically valid over all possible ways one could define a proof-theoretic system.

Tehcinically we can write this as:
\begin{defn}\edit{
$\varphi$ is a \emph{generalised proof-theoretically valid formula} if for every proof-theoretic system $\mathfrak{S}$ and $S\in\mathfrak{S},$ it follows that there is an $S$-valid proof of $\varphi$ where all acceptable extensions of $S$ are in $\mathfrak{S}$.}
\end{defn}

Then it can be proven that:
\begin{thm}[Revised Prawitz's Conjecture ([citation removed to preserve blind review])]
For all $\varphi$ in the language of propositional logic, $\varphi$ is a generalised proof-theoretically valid formula $\Leftrightarrow$ $\varphi$ is an intuitionistic validity.
\end{thm}

The modified Prawitz’s conjecture is then that this definition is equivalent to intuitionistic logic.  And the result mentioned above (that there is a counter example to every proof-theoretic validity) is sufficient to prove the correctness of the conjecture. }

This is in contrast with how proof-theoretic validity was treated in \textcite{Piecha2019-sq} where the sets of atomic rules, not the proof-theoretic systems, are treated analogously with models.  This proposal avoids the problem we have identified above and is not subject to \textcite{Piecha2019-sq} general result about superintuitionistic validates. \textcite[\S 4]{Piecha2019-sq} are aware that it is possible to generalises definitions as discussed here, as shown in their discussion of \textcite{Goldfarb2016-tk}'s approach which can be seen as a restricted version of the approach proposed here. However, what we hope to have shown here is that there are reasons to prefer the more general definition of proof-theoretic consequence over the more restricted and superintuitionistic notions. 

\section{Conclusion}
Let us reiterate the argument of the paper. Proof-theoretic validity is supposed to be a formal method for finding those inferences that follow from intuitionistic introduction rules. It aims to explicate Gentzen’s claim that intuitionistic elimination rules are consequences of intuitionistic introduction rules, and it was conjectured by Prawitz and Dummett that no intuitionistically invalid inferences are proof-theoretically valid. Instead, it has been shown that superintuitionistic inferences are proof-theoretically valid. This may look like a bad situation for proof-theoretic validity because it seems tenuous to claim that these inferences follow from the introduction rules.

In the end, however, a careful examination of the definition of proof-theoretic validity and the proof of superintuitionistic validities highlighted the vital role of the treatment of the atomic formulas. Whether an instance of a superintuitionistic rule was valid or not turned out to depend on whether there was an atomic rule corresponding to the formulas it contained. On this basis, it was argued that it is the treatment of atomic formulas that is leading to the validity of superintuitionistic inferences.  And we propose a revised Prawitz's conjecture that is true for proof-theoretic validity defined over all proof-theoretic systems rather than one.


\printbibliography
\end{document}